\documentstyle[twoside,psfig,lscape]{article}

\oddsidemargin=\evensidemargin 
\catcode`\@=11 \long\def\@makefntext#1{ \protect\noindent \hbox to
3.2pt {\hskip-.9pt
$^{{\eightrm\@thefnmark}}$\hfil}#1\hfill}       

\def\ps@myheadings{\let\@mkboth\@gobbletwo      
\def\@oddhead{\hbox{}
\rightmark\hfil\eightrm\thepage}
\def\@oddfoot{}\def\@evenhead{\eightrm\thepage\hfil
\leftmark\hbox{}}\def\@evenfoot{}
\def\sectionmark##1{}\def\subsectionmark##1{}}

\catcode`@=11                        
\def\ps@plain{\let\@mkboth\@gobbletwo
     \def\@oddhead{}\def\@oddfoot{\eightrm\hfil\thepage
     \hfil}\def\@evenhead{}\let\@evenfoot\@oddfoot}

\newcounter{sectionc}\newcounter{subsectionc}\newcounter{subsubsectionc}
\renewcommand{\section}[1] {\vspace{12pt}\addtocounter{sectionc}{1}
\setcounter{subsectionc}{0}\setcounter{subsubsectionc}{0}\noindent
    {\tenbf\thesectionc. #1}\par\vspace{5pt}}
\renewcommand{\subsection}[1] {\vspace{12pt}\addtocounter{subsectionc}{1}
    \setcounter{subsubsectionc}{0}\noindent
    {\bf\thesectionc.\thesubsectionc.
    {\kern1pt \bfit #1}}\par\vspace{5pt}}
\renewcommand{\subsubsection}[1] {\vspace{12pt}
    \addtocounter{subsubsectionc}{1}
    \noindent
    {\tenrm\thesectionc.\thesubsectionc.\thesubsubsectionc. {\kern1pt
    \it #1}}\par\vspace{5pt}}

\newcounter{appendixc}
\newcounter{subappendixc}[appendixc]
\newcounter{subsubappendixc}[subappendixc]

\renewcommand{\appendix}[1] {\vspace{12pt}  
    \refstepcounter{appendixc}      
    \setcounter{figure}{0}
    \setcounter{table}{0}
    \setcounter{lemma}{0}
    \setcounter{theorem}{0}
    \setcounter{corollary}{0}
    \setcounter{definition}{0}
    \setcounter{equation}{0}
    \renewcommand{\thefigure}{\Alph{appendixc}.\arabic{figure}}
    \renewcommand{\thetable}{\Alph{appendixc}.\arabic{table}}
    \renewcommand{\theappendixc}{\Alph{appendixc}}
    \renewcommand{\thelemma}{\Alph{appendixc}.\arabic{lemma}}
    \renewcommand{\thetheorem}{\Alph{appendixc}.\arabic{theorem}}
    \renewcommand{\thedefinition}{\Alph{appendixc}.\arabic{definition}}
    \renewcommand{\thecorollary}{\Alph{appendixc}.\arabic{corollary}}
    \renewcommand{\theequation}{\Alph{appendixc}.\arabic{equation}}
    \noindent{\tenbf Appendix \theappendixc #1}\par\vspace{5pt}}

\newcommand{\smalllineskip}{\baselineskip=10pt}

\newcommand{\copyrightheading}[1]
    {\vspace*{-2.5cm}\smalllineskip{\flushleft
    {\footnotesize }\\
    {\footnotesize \copyright\kern2pt }\\
         }}

\def\keywords#1{{
    \centering{\begin{minipage}{4.5in}\footnotesize\baselineskip=10pt
    {\footnotesize\it Keywords}\/: #1
    \end{minipage}}\par}}

\renewenvironment{thebibliography}[1]
    {\frenchspacing
     \ninerm\baselineskip=11pt
     \begin{list}{[\arabic{enumi}]}
    {\usecounter{enumi}\setlength{\parsep}{0pt}
     \setlength{\leftmargin 13.7pt}{\rightmargin 0pt} 
     \setlength{\itemsep}{0pt} \settowidth
    {\labelwidth}{[#1]}\sloppy}}{\end{list}}

\newcounter{itemlistc}
\newcounter{romanlistc}
\newcounter{alphlistc}
\newcounter{arabiclistc}

\newcommand{\fcaption}[1]{
        \refstepcounter{figure}
        \setbox\@tempboxa = \hbox{\footnotesize Fig.~\thefigure. #1}
        \ifdim \wd\@tempboxa > 5in
           {\begin{center}
        \parbox{5in}{\footnotesize\smalllineskip Fig.~\thefigure. #1}
            \end{center}}
        \else
             {\begin{center}
             {\footnotesize Fig.~\thefigure. #1}
              \end{center}}
        \fi}

\newcommand{\tcaption}[1]{
        \refstepcounter{table}
        \setbox\@tempboxa = \hbox{\footnotesize Table~\thetable. #1}
        \ifdim \wd\@tempboxa > 5in
           {\begin{center}
        \parbox{5in}{\footnotesize\smalllineskip Table~\thetable. #1}
            \end{center}}
        \else
             {\begin{center}
             {\footnotesize Table~\thetable. #1}
              \end{center}}
        \fi}

\def\pmb#1{\setbox0=\hbox{#1}
    \kern-.025em\copy0\kern-\wd0
    \kern.05em\copy0\kern-\wd0
    \kern-.025em\raise.0433em\box0}

\def\fnt#1#2{\footnotetext{\kern-.3em
    {$^{\mbox{\scriptsize #1}}$}{#2}}}

\font\tenrm=cmr10  \font\tenbf=cmbx10
\font\bfit=cmbxti10 at 10pt \font\ninerm=cmr9 
 \font\eightrm=cmr8

\newtheorem{definition}{Definition}

\def\@begintheorem#1#2{\trivlist    
    \item[\hskip\labelsep{\bf #1\ #2.}]}
\def\@opargbegintheorem#1#2#3{\trivlist
    \item[\hskip\labelsep{\bf #1\ #2\ (#3).}]}

        {\setcounter{itemlistc}{0}      
     \begin{list}{$\bullet$}        
    {\usecounter{itemlistc}         
     \leftmargin10pt           
     \setlength{\parsep}{0pt}
     \setlength{\itemsep}{0pt}     
    }}{\end{list}}

    {\setcounter{romanlistc}{0}     
     \begin{list}{$($\roman{romanlistc}$)$} 
    {\usecounter{romanlistc}        
     \leftmargin18pt 
     \setlength{\parsep}{0pt}
     \setlength{\itemsep}{0pt}  
     \settowidth{\labelwidth}{#1}
    }}{\end{list}}

    {\setcounter{enumii}{0}         
     \begin{list}{$($\alph{enumii}$)$}  
    {\usecounter{enumii}            
     \leftmargin18pt        
     \setlength{\parsep}{0pt}
     \setlength{\itemsep}{0pt}  
     \settowidth{\labelwidth}{#1}
    }}{\end{list}}

\def\qed{\hbox{${\vcenter{\vbox{            
   \hrule height 0.4pt\hbox{\vrule width 0.4pt height 6pt
   \kern5pt\vrule width 0.4pt}\hrule height 0.4pt}}}$}}

\def\theequation{\thesectionc.\arabic{equation}}  

\begin{document}

\markboth{Louis H. Kauffman, Slavik V. Jablan,  Ljiljana Radovi\' c, Radmila Sazdanovi\'
c} {Reduced relative Tutte, Kauffman bracket and Jones polynomials of virtual link families}

\centerline{\bf REDUCED RELATIVE TUTTE, KAUFFMAN
BRACKET AND}
\centerline{\bf JONES POLYNOMIALS OF VIRTUAL LINK FAMILIES}
\bigskip

\centerline{\footnotesize LOUIS H. KAUFFMAN, SLAVIK V. JABLAN$^\dag $,}
\centerline{\footnotesize LJILJANA RADOVI\' C$^{\dag
\dag }$, RADMILA SAZDANOVI\' C$^{\dag \dag \dag}$}

\bigskip

\centerline{\footnotesize\it University of Illinois at Chicago,
Department of Mathematics,}\centerline{\footnotesize\it Statistics
and Computer Science (m/c 249),}\centerline{\footnotesize\it 851
South Morgan Street, Chicago,}\centerline{\footnotesize\it Illinois
60607-7045, USA} \centerline{\footnotesize\it kauffman@uic.edu}

\bigskip

\centerline{\footnotesize\it The Mathematical Institute$^\dag $, Knez
Mihailova 36,}\centerline{\footnotesize\it P.O.Box 367, 11001
Belgrade, Serbia} \centerline{\footnotesize\it sjablan@gmail.com}

\bigskip

\centerline{\footnotesize\it Faculty of Mechanical
Engineering$^{\dag \dag}$, A.~Medvedeva 14,}
\centerline{\footnotesize\it 18 000 Ni\v s, Serbia} \centerline{\footnotesize\it ljradovic@gmail.com}

\bigskip

\centerline{\footnotesize\it University of Pennsylvania$^{\dag \dag \dag}$, 209 South 33rd Street,}
\centerline{\footnotesize\it Philadelphia, PA 19104-6395, USA} \centerline{\footnotesize\it radmilas@gmail.com}

\bigskip

\begin{abstract}
\noindent This paper contains general formulae for the reduced
relative Tutte, Kauffman bracket and Jones polynomials of families
of virtual knots and links given in Conway notation and discussion of a counterexample to
 the $Z$-move conjecture of Fenn, Kauffman and Manturov.
\end{abstract}

\bigskip

\keywords{Conway notation; virtual knot; Tutte polynomial; Bollob\'
as-Riordan polynomial; relative Tutte polynomial; reduced relative
Tutte polynomial; Jones polynomial; Kauffman bracket polynomial;
link family.}

\section{Introduction}

According to Thistlethwaite's Theorem \cite{1} and Kauffman's extension of it \cite{2}, the Jones polynomial and Kauffman bracket polynomial of any classical  knot or link (shortly $KL$) can be computed from the
(signed) graph of the $KL$ \cite{2} via the Tutte polynomial for signed graphs.
Without computing Tutte polynomials, the Jones polynomials of the
diagrams with at most 5 twist regions in which crossings can be
inserted are computed by D. Silver, A. Stoimenow, and S. Williams
\cite{3}, and the general method for computing Kauffman bracket
polynomial of such diagrams is outlined by A. Champanerkar and I.
Kofman \cite{4}.

The general formulae for Tutte and Jones polynomials for families of
classical $KL$s with at most 5 twisting regions, given in Conway
notation, are derived in \cite{5}.

In Subsections 2.1-2.3 we compute Kauffman bracket and Jones polynomials for several families of virtual knots using the Tutte polynomials of their
signed graph. In Section 4 we use Vassily Manturov's parity bracket to find a new
counterexample to the $Z$-move conjecture, and in Section 5 plot zeros of Jones polynomials
of certain families of virtual knots and links.

Knots and links (or shortly $KL$s) can be given in Conway notation
\cite{6,7,8,9}. For readers unfamiliar with it, we explain Conway notation, introduced in Conway's seminal paper
\cite{6} published in 1967, and effectively used since (e.g.,
\cite{7,8,9}). Conway symbols of knots with up to $10$ crossings and
links with at most $9$ crossings are given in the Appendix of the
book \cite{7}.

\begin{figure}[th]
\centerline{\psfig{file=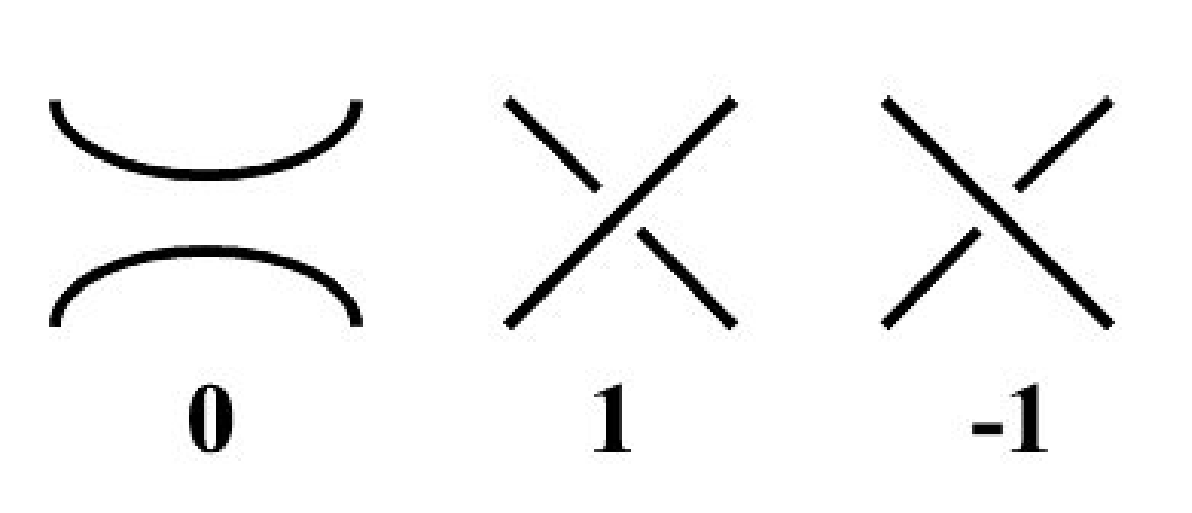,width=2.00in}} \vspace*{8pt}
\caption{The elementary tangles. \label{f1}}
\end{figure}

\begin{figure}[th]
\centerline{\psfig{file=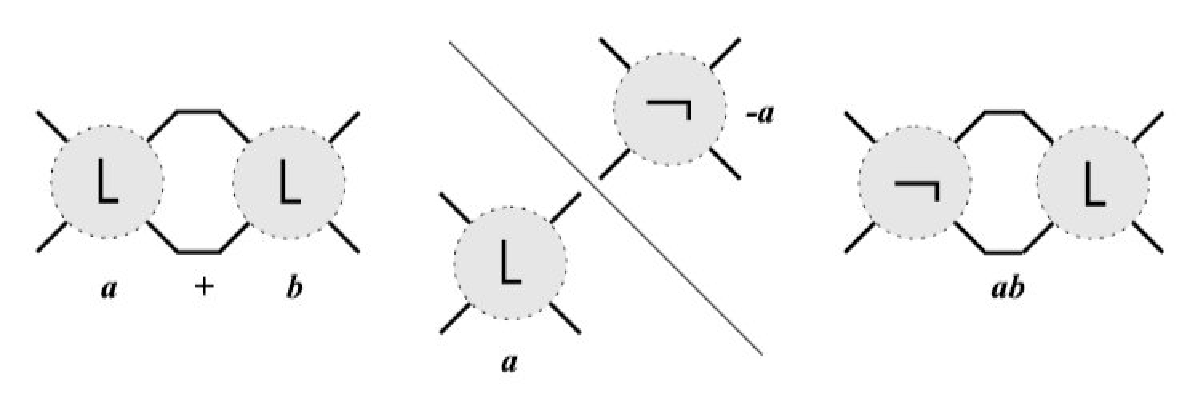,width=3.80in}} \vspace*{8pt}
\caption{A sum and product of tangles. \label{f2}}
\end{figure}

\begin{figure}[th]
\centerline{\psfig{file=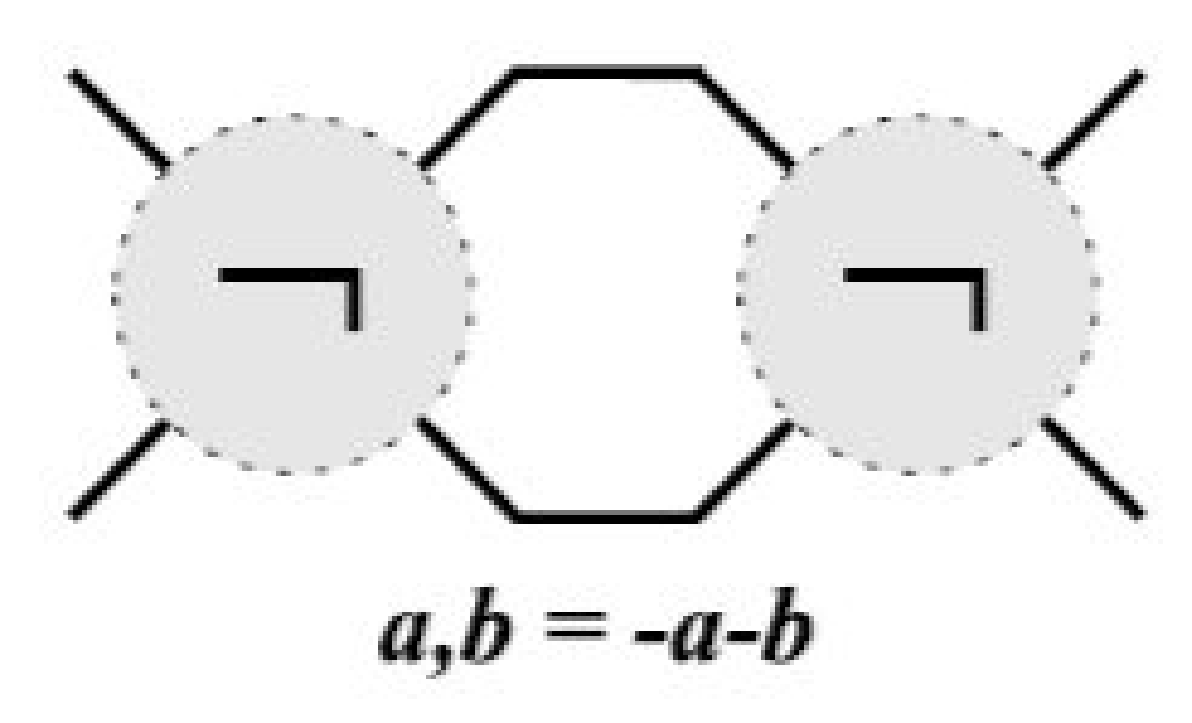,width=1.40in}} \vspace*{8pt}
\caption{Ramification of tangles. \label{f3}}
\end{figure}

The main building blocks in Conway notation are elementary
tangles. We distinguish three elementary tangles,  shown in Fig.
\ref{f1} and denoted by 0, 1 and $-1$. All other tangles can be
obtained by combining elementary tangles, while 0 and 1 are
sufficient for generating alternating knots and links.
Elementary tangles can be combined by following operations: {\it
sum}, {\it product}, and {\it ramification} (Figs.
\ref{f2}-\ref{f3}). Given tangles $a$ and $b$, image of $a$ under
reflection with mirror line NW-SE is denoted by $-a$, and sum is
denoted by $a+b$. The product $a\,b$ is defined as $a\,b = -a+b$, and
{\it ramification} by $(a,b) = -a-b$.

\begin{figure}[th]
\centerline{\psfig{file=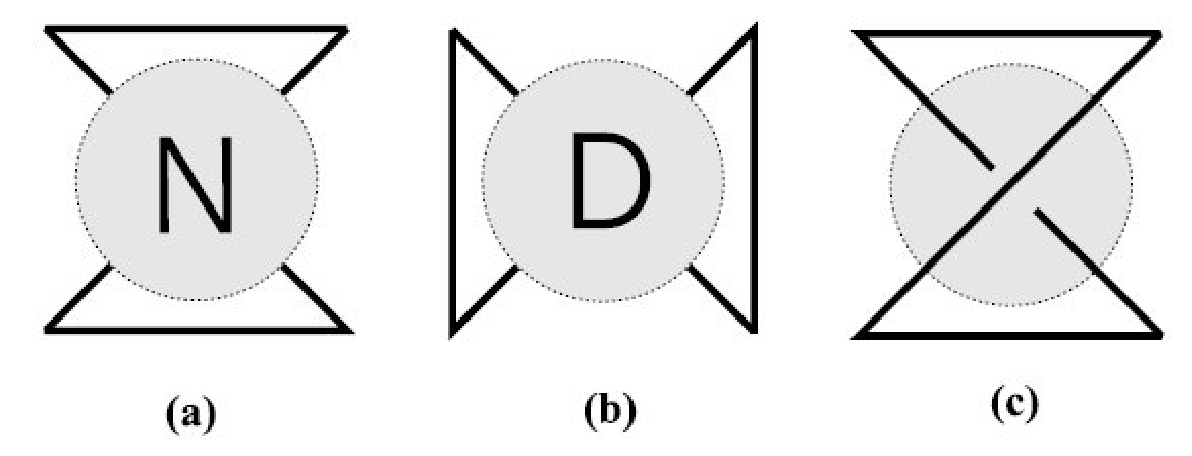,width=2.40in}} \vspace*{8pt}
\caption{(a) Numerator closure; (b) denominator closure; (c) basic
polyhedron $1^*$. \label{f5}}
\end{figure}

The tangle can be closed in two ways (without
introducing additional crossings): by joining in pairs NE and NW, and
SE and SW ends of a tangle to obtain a {\it numerator closure};
or by joining in pairs NE and SE, and NW and SW ends we obtain a {\it
denominator closure} (Fig. \ref{f5}a,b).

\begin{definition}
A {\it rational tangle} is any finite product of elementary tangles.
A {\it rational $KL$} is a numerator closure of a rational tangle.
\end{definition}

\begin{definition}
A tangle is {\it algebraic} if it can be obtained from elementary
tangles using the operations of sum and product. A $KL$ is {\it
algebraic} if it is a numerator closure of an algebraic tangle.
\end{definition}

\begin{definition}
For a link or knot $L$ given in an unreduced\footnote{The Conway
notation is called unreduced if in symbols of nonalgebraic links
elementary tangles 1 in single vertices are not omitted.} Conway
notation $C(L)$ denote by $S$ a set of numbers in the Conway symbol
excluding numbers denoting basic polyhedron and zeros (determining
the position of tangles in the vertices of polyhedron) and let
$\tilde S=\{a_1,a_2, \ldots, a_k\}$ be a non-empty subset of $S$.
Family $F_{\tilde S}(L)$ of knots or links derived from $L$ consists
of all knots or links $L'$ whose Conway symbol is obtained by
substituting all $ a_i\neq \pm 1$, by $sgn(a_i) |a_i+k_{a_i}|$,
$|a_i+k_{a_i}| >1$, $k_{a_i} \in Z$.
\end{definition}

An infinite subset of a family is called a {\it subfamily}. If all
$k_{a_i}$ are even integers, the number of components is preserved
within the corresponding subfamilies, i.e., adding full-twists
preserves the number of components inside the subfamilies.

\begin{definition}
A link given by Conway symbol containing only tangles $\pm 1$ and  $\pm 2$
is called a {\it source link}. A link given by Conway
symbol containing only tangles $\pm 1$, $\pm 2$, or $ \pm 3$ is
called a {\it generating link}.
\end{definition}

\begin{figure}[th]
\centerline{\psfig{file=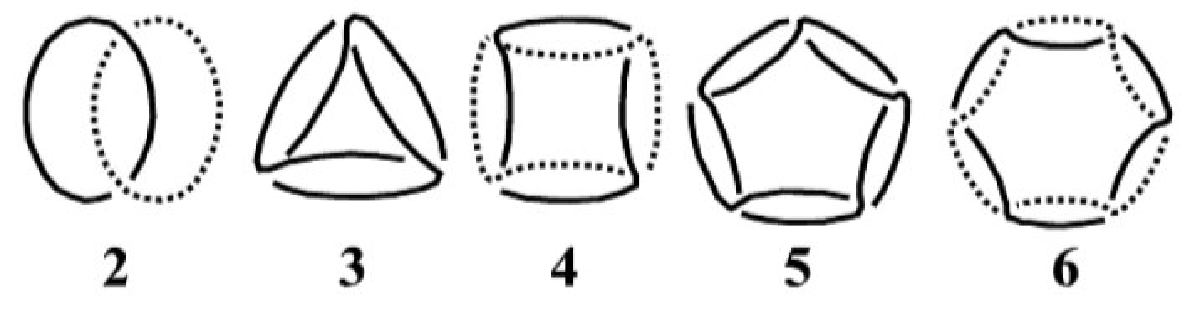,width=3.00in}} \vspace*{8pt}
\caption{Hopf link $2_1^2=2$ and its family $p$ ($p=2,3,\ldots $).
\label{f12}}
\end{figure}

For example, Hopf link $2$ (link $2_1^2$ in Rolfsen's notation) is
the source link of the simplest link family $p$ ($p=2,3,\ldots $)
(Fig. \ref{f12}), and Hopf link and trefoil $3$ (knot $3_1$ in the
classical notation) are generating links of this family. A family of
$KL$s is usually derived from its source link by substituting $a_i
\in \tilde S$, $a_i=\pm 2$, by $sgn(a_i) (2+k)$, $k=1,2,3,\ldots $
(see Def. 1.5 and Def. 1.6).

Analogous to the Conway notation for classical $KL$s
we use extended Conway notation for virtual $KL$s, by adding to the list of the elementary tangles $0$, $1$, $-1$ the elementary tangle $i$ denoting a virtual crossing. The extended
Conway notation for virtual links differs from the standard one in
the following way:
\begin{itemize}
  \item virtual crossings are denoted by $i$.
  \item sequence of $n$ classical crossings $1,\ldots ,1$ (positive $n$-twist) is denoted by $1^n$
  \item sequence of $n$ classical negative crossings  $-1,\ldots ,-1$ (negative $n$-twist) is denoted by
$(-1)^n$
\end{itemize}

The convention introduced in the Conway notation extended to virtual links is clear from the fact that every  positive $n$-twist can be denoted as $1,\ldots ,1=1^n$, and  every negative $n$-twist by $-1,\ldots ,-1.$ This convention extends to virtual knot families. For example, the simplest $KL$ family $2_1^2$, $3_1$, $4_1^2$, $5_1$, $\ldots$ of knots and links, consisting from Hopf link $2$, trefoil $3$, $\ldots $ can be denoted by general symbol $p$ ($p\ge 2$). From them, by substituting in each of them one classical crossing by virtual one, we obtain a family of virtual $KL$s $(i,1)$, $(i,1,1=i,1^2)$, $(i,1,1,1=i,1^3)$, $(i,1,1,1,1=i,1^4)$, $\ldots$, and in general $i,1^{p-1}$ ($p\ge 2$).

In comparison with other $KL$ notations (Gauss codes, Kamada codes, PD-notation, {\it etc}.), Conway notation extended to virtual links is much shorter than any other notation for virtual $KL$s and the
most suitable for utilizing the notion of families of knots and
links and analyzing how knot and link properties change inside
families. E.g, the Gauss code of the trefoil with one virtual crossing is $U1+O2+O1+U2+$, its Kamada code is $\{\{\{1,3\},\{2,0\}\},\{1,1\}\}$, its PD-code is PD[X[3,1,4,2],X[2,4,3,1]], and its Conway symbol is $(1,1,i)=(1^2,i)$, meaning that in the classical trefoil knot $3=(1,1,1)=1^3$ the last classical crossing  $1$ is substituted by the virtual crossing $i$.

The relative Tutte polynomial of colored graphs is introduced in the
paper \cite{10} by Y.~Diao and G.~Hetei. An alternative approach to
the computation of Bollob\' as-Riordan polynomial and Kauffman
bracket polynomial of virtual $KL$s, via ribbon graphs, is given by
S. Chmutov and I. Pak \cite{11}, and S. Chmutov \cite{12}.

Zeroes of the Jones polynomials of different knots and their plots are computed
in \cite{13,14,15,16}. Zeroes of Jones polynomials corresponding to
different families of classical (alternating and non-alternating)
$KL$s, called ``portraits of the families'', are given in \cite{5}.

First we restrict our attention to graphs corresponding only to
alternating virtual $KL$s, hence we consider graphs corresponding to
virtual $KL$s with all edges labeled by $+$ or 0 (see \cite{10},
Sections 2, 5), with variables $X_+=X$, $Y_+=Y$, $x_+=x$, $y_+=y$, $X_0=1$, $Y_0=1$, $x_0=x$, $y_0=y$. In
this setting, we have the following recursive formula for computing
the relative Tutte polynomial:

$$T(G)=\cases {  $if$ \quad  G \quad  $is a graph with no edges$ \quad   \quad   \quad  \quad   \quad   \quad  \quad  \quad  \quad  \quad \quad  \quad \quad  \quad  (1)  \cr
X_{\lambda }T(G/e) \quad $if$ \quad   e \quad  $is a bridge$ \quad   \quad \quad     \quad  \quad   \quad   \quad  \quad  \quad  \quad  \quad \quad  \quad \quad  \quad  (2)  \cr
Y_{\lambda }T(G-e) \quad $if$ \quad  e \quad $is a loop$ \quad   \quad   \quad  \quad   \quad   \quad  \quad  \quad  \quad  \quad \quad  \quad \quad  \quad \quad $\,\,$ (3)  \cr
y_{\lambda }T(G-e)+x_{\lambda }T(G/e) \quad $if$ \quad  e  \quad $is neither a loop or a bridge$   $\,\,\,$ \quad  (4)  \cr
                }$$

\noindent where $G-e$ is the graph obtained from $G$ by deleting the
edge $e$, and $G/e$ is the graph obtained from $G$ by contracting
$e$, and $\lambda $ is the color of $e$.

Notice that formulae for the relative Tutte polynomial of
non-alternating virtual $KL$s, can be obtained from general formulae
for alternating virtual $KL$s by substituting negative values of
parameters. For a graph $G$ and its dual $G'$ variables in their corresponding Tutte polynomials $T(G)$ and $T(G')$ change their places.
For graphs $G$ and $H$, $T(G\ast H)=T(G)T(H)$, where $G\ast H$ is the join of $G$ and $H$ at a vertex.

According to the Theorem 5.4 \cite{10}, the relative Tutte
polynomial $T_H(G)$ (see Section 5.1 \cite{10}) corresponding to a
virtual link diagram $K$ gives the Kauffman bracket through the
following variable substitution:

$$X\rightarrow -A^{-3}, \quad Y\rightarrow -A^3, \quad x\rightarrow A,
\quad y\rightarrow A^{-1}, \quad d\rightarrow -(A^2+A^{-2}),$$

\noindent and the Jones polynomial of $K$ is obtained from the
Kauffman bracket polynomial by substituting $A=t^{-{1\over 4}}$.

\begin{figure}[th]
\centerline{\psfig{file=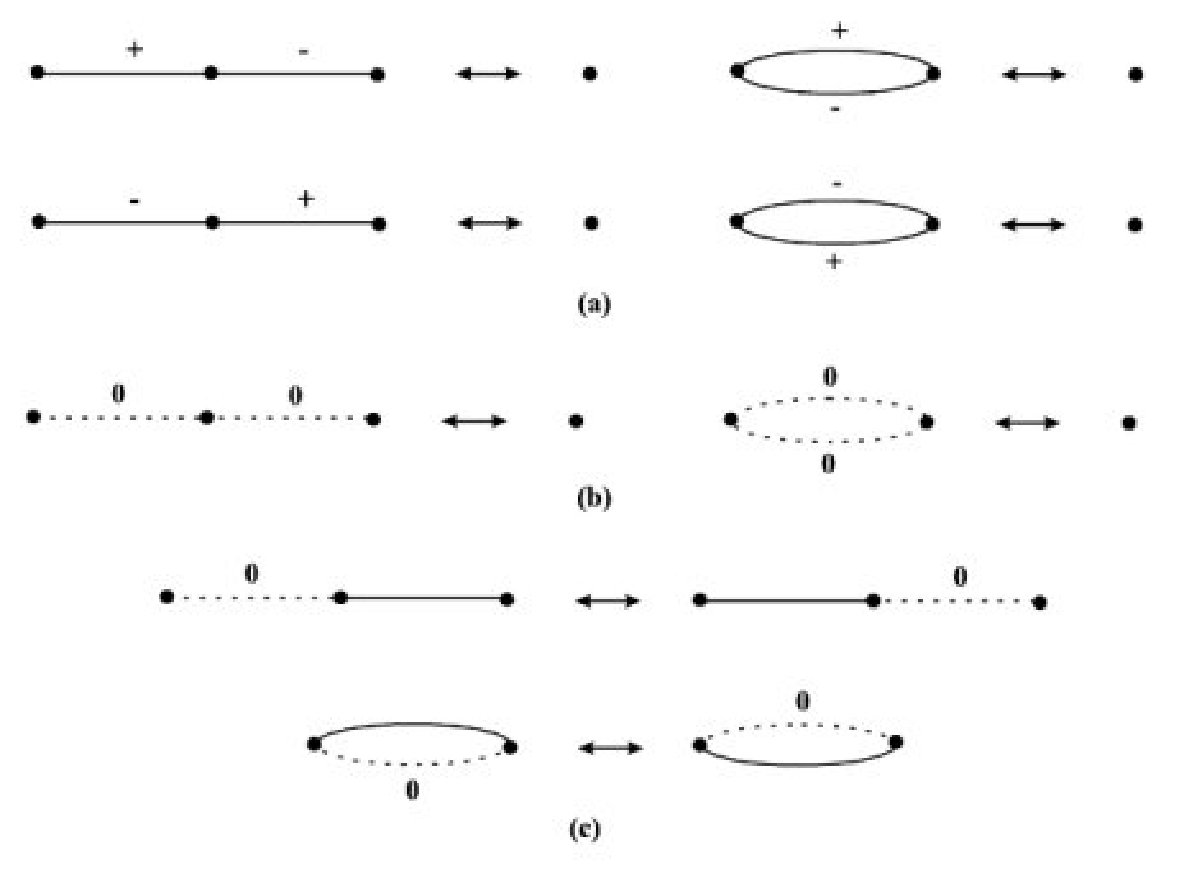,width=3.40in}} \vspace*{8pt}
\caption{(a) Reidemeister move II for classical crossings; (b)
Reidemeister move II for virtual crossings; (c) $Z$-move (virtualization) on $KL$
graphs. Non-signed edges can have sign $+$ or $-$. \label{f1.1}}
\end{figure}

\section{Reduction rules}
\begin{figure}[th]
\centerline{\psfig{file=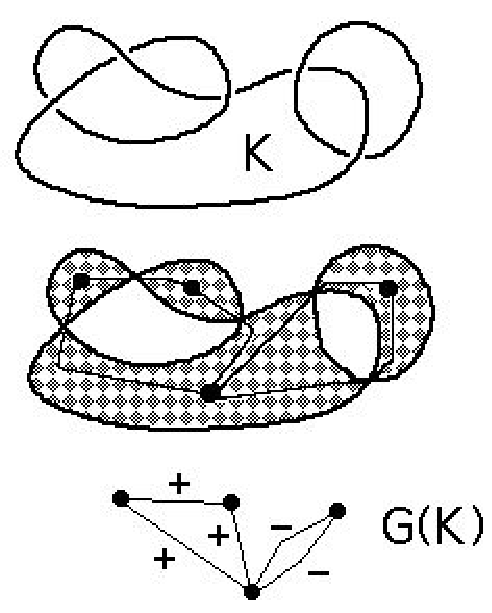,width=2.00in}} \vspace*{8pt}
\caption{Signed Graph of a Knot. \label{checker}}
\end{figure}

\begin{figure}[th]
\centerline{\psfig{file=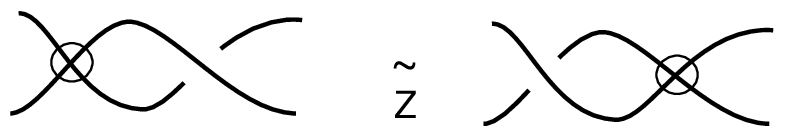,width=2.00in}} \vspace*{8pt}
\caption{The $Z$ - Move. \label{zmove}}
\end{figure}

In this section we consider certain reduction and transformation rules relating graphs,
knots and links. For this we first recall the construction of the signed graph of a knot or link diagram.
View Fig. \ref{checker}. In this figure we illustrate how one first forms the two-colored checkerboard for the knot or link diagram, and then associates a signed graph $G(K)$ to this checkerboard by assigning a node to each shaded region of the checkerboard and and edge to each crossing that is incident to two regions. We label the edges of the resulting graph with $+1$, $-1$ and $0$ according as the crossing is classical and positive in relation to the edge, classical and negative in relation to the edge, and virtual in relation to the edge. Virtual crossings are neither over nor under and they are indicated on the knot diagram by a flat crossing with a circle around it as in Fig. \ref{zmove}. This figure also illustrates the
{\em $Z$-move} (or {\em virtualization},  a move on virtual knot and link diagrams that does not effect the evaluation of the
Kauffman bracket or the Jones polynomial. The reader who would like more information about virtual
knot theory can consult \cite{IVKT}. In Fig. \ref{f1.1}c we illustrate the effect of the $Z$-Move on
the graph of a virtual knot or link. Note that we have indicated virtual crossings by corresponding graph edges that are labeled with a $0$ and are dotted edges. We have not indicated any signs in this
figure.

Since the Kauffman bracket and Jones polynomial of a virtual link $L$ are
invariant under classical and virtual Reidemeister move II and $Z$-move
(\cite{17}, Lemma 7) on graphs, we can simplify graphs before
computations. Graphs of virtual $KL$s can be simplified using the
following reduction rules:

\begin{enumerate}
\item Reidemeister move II for classical crossings (Fig. \ref{f1.1}a);
\item Reidemeister move II for virtual crossings (Fig. \ref{f1.1}b);
\item $Z$-move (Fig. \ref{f1.1}c) on $KL$ graphs.
\end{enumerate}

\noindent where $0$-edges are denoted by broken lines. The relative
Tutte polynomial obtained after this reduction will be called {\it
reduced relative Tutte polynomial}.

A graph is called completely reduced if none of the reduction rules can be applied further. Twists obtained by reductions can contain crossings of the same sign, and at most one virtual crossing, so their corresponding parts of reduced graphs can contain only sequences of edges of the same sign or sets of multiple edges of the same sign, and at most one 0-edge each.

\begin{figure}[th]
\centerline{\psfig{file=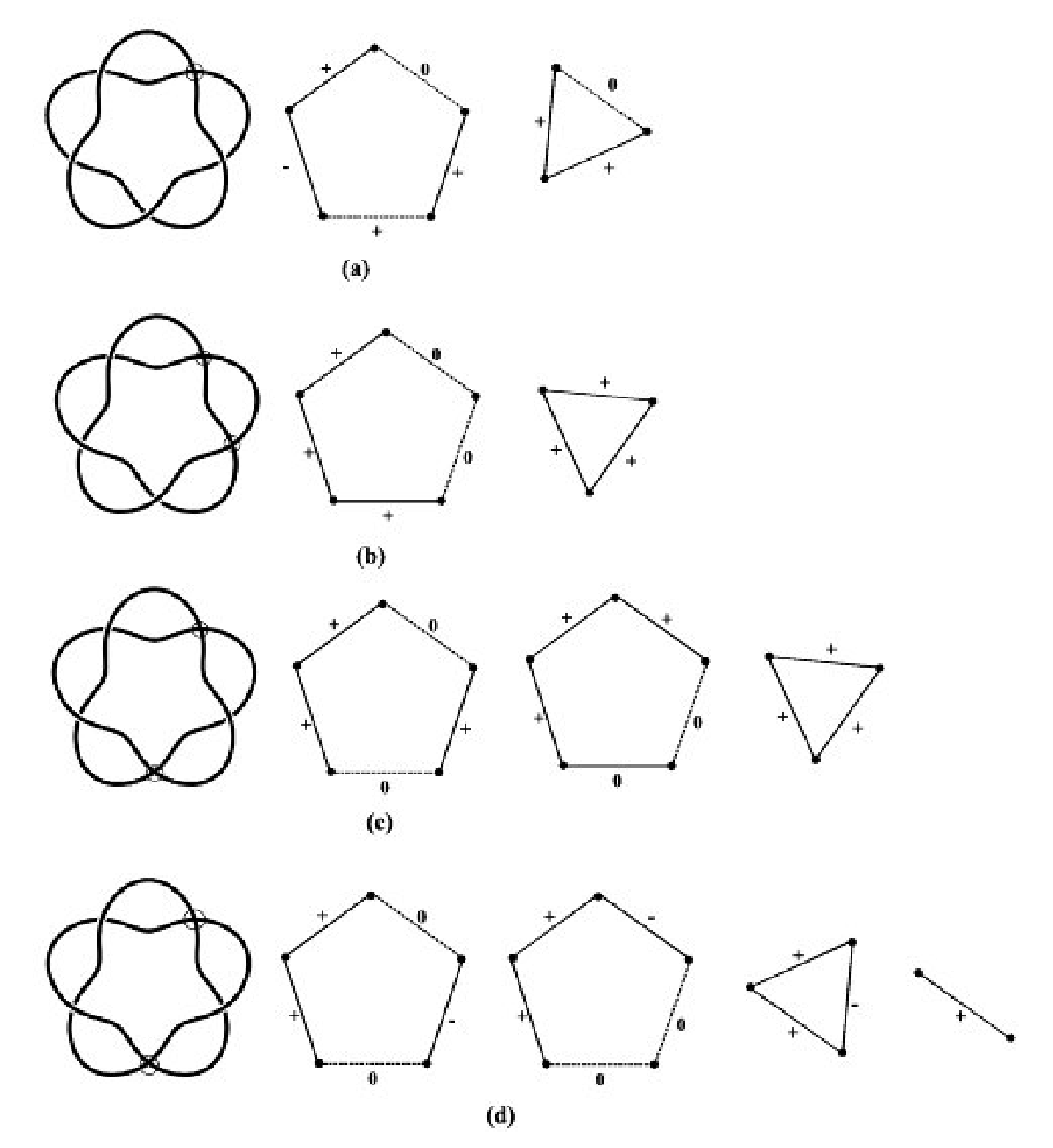,width=4.20in}} \vspace*{8pt}
\caption{Reduction of the graphs of knots (a) $1,-1,1,1,i$; (b)
$1,1,1,i,i$; (c)  $1,1,i,1,i$. (d)  $1,1,i,-1,i$. \label{reduct}}
\end{figure}

As a simple example of the application of reduction rules to graphs of
virtual $KL$s consider virtual knots $(1,-1,1,1,i)$, $(1,1,1,i,i)$, and
$(1,1,i,1,i)$ (Fig. \ref{reduct}). From the graph of the first knot, by
the Reidemeister move II for classical crossings we conclude that its
graph reduces to the graph of a trefoil with one virtual crossing (Fig.
\ref{reduct}a). From the graph of the second knot, by the Reidemeister
move II for virtual crossings, we obtain the graph of the trefoil (Fig.
\ref{reduct}b). From the graph of the third knot, by $Z$-move and
Reidemeister move II for virtual crossings, we conclude that its Kauffman
bracket and Jones polynomial are same as those of the trefoil $1,1,1$
($1,1,i,1,i\rightarrow 1,1,1,i,i\rightarrow 1,1,1=3$) (Fig.
\ref{reduct}c). Using graph reduction, $Z$-move followed by Reidemeister
move II for virtual crossings and Reidemeister move II for classical
crossings,  we conclude that the Kauffman bracket and Jones polynomial of
the knot $1,1,i,-1,i$ are unit ($1,1,i,-1,i\rightarrow
1,1,-1,i,i\rightarrow 1,1,-1=1$)(Fig. \ref{reduct}d).  A graph is called
completely reduced if none of the reduction rules can be applied
further.

\subsection{The family \it{\textbf{p}}}

Consider the family \emph{\textbf{p}} ($p\ge 2$), which
consists of classical knots and and links of the form $1^p$. Graphs
corresponding to links of this family are cycles of length $p$,
satisfying the following recursion:

$$T(G(1^p))=yX^{p-1}+xT(G(1^{p-1}))$$

\noindent with $T(G(1^2))= yX+xY$, so the general formula for the
reduced relative Tutte polynomial of the graph $G((1^p))$ is

$$T(G(1^p))={(x^p-X^p) \over {x-X}}y+(Y-y)x^{p-1}.$$

\begin{figure}[th]
\centerline{\psfig{file=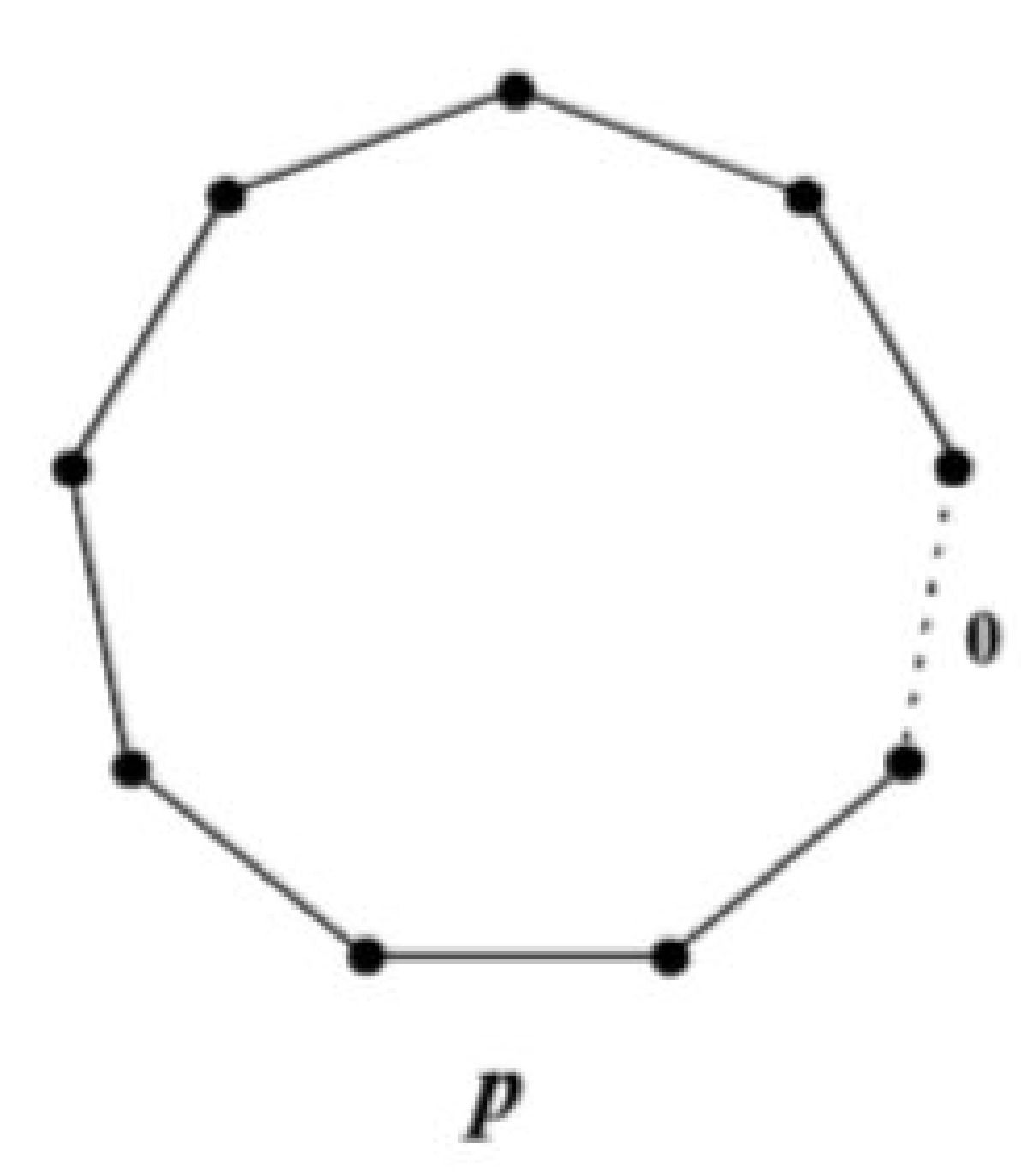,width=1.00in}} \vspace*{8pt}
\caption{Graph $G(i,1^p)$.\label{f1.2}}
\end{figure}

Virtual $KL$s of the form $(i,1^p)$ which belong to the same family, correspond to
the reduced graphs shown on Fig. \ref{f1.2}, and their Tutte polynomials satisfy
 the following recursion:

$$T(G(i,1^p))=yX^{p-1}+xT(G(i,1^{p-1}))$$

\noindent with $T(G(i,1))= x+y$,  so the general formula for the
reduced relative Tutte polynomial of the graph $G((i,1^p))$ is

$$T(G(i,1^p))=x^p+{(x^p-X^p) \over {x-X}}y.$$

As a corollary of this general formula we obtain reduced relative Tutte
polynomials for positive or negative values of the parameter $p$,
and Kauffman bracket and Jones polynomials. For example, for $p=2$
we obtain the reduced relative Tutte polynomial $T(G(i,1^2))=x^2 +
xy + Xy$ of the virtual trefoil $(i,1,1)$, and for $p=-2$ the
reduced relative Tutte polynomial $T(G(i,(-1)^2))={1\over x^2} -
{y\over xX^2} - {y\over x^2X}$ of its mirror image $(i,-1,-1)$.

\subsection{The family \it{\textbf{p\,q}}}

\begin{figure}[h]
\centerline{\psfig{file=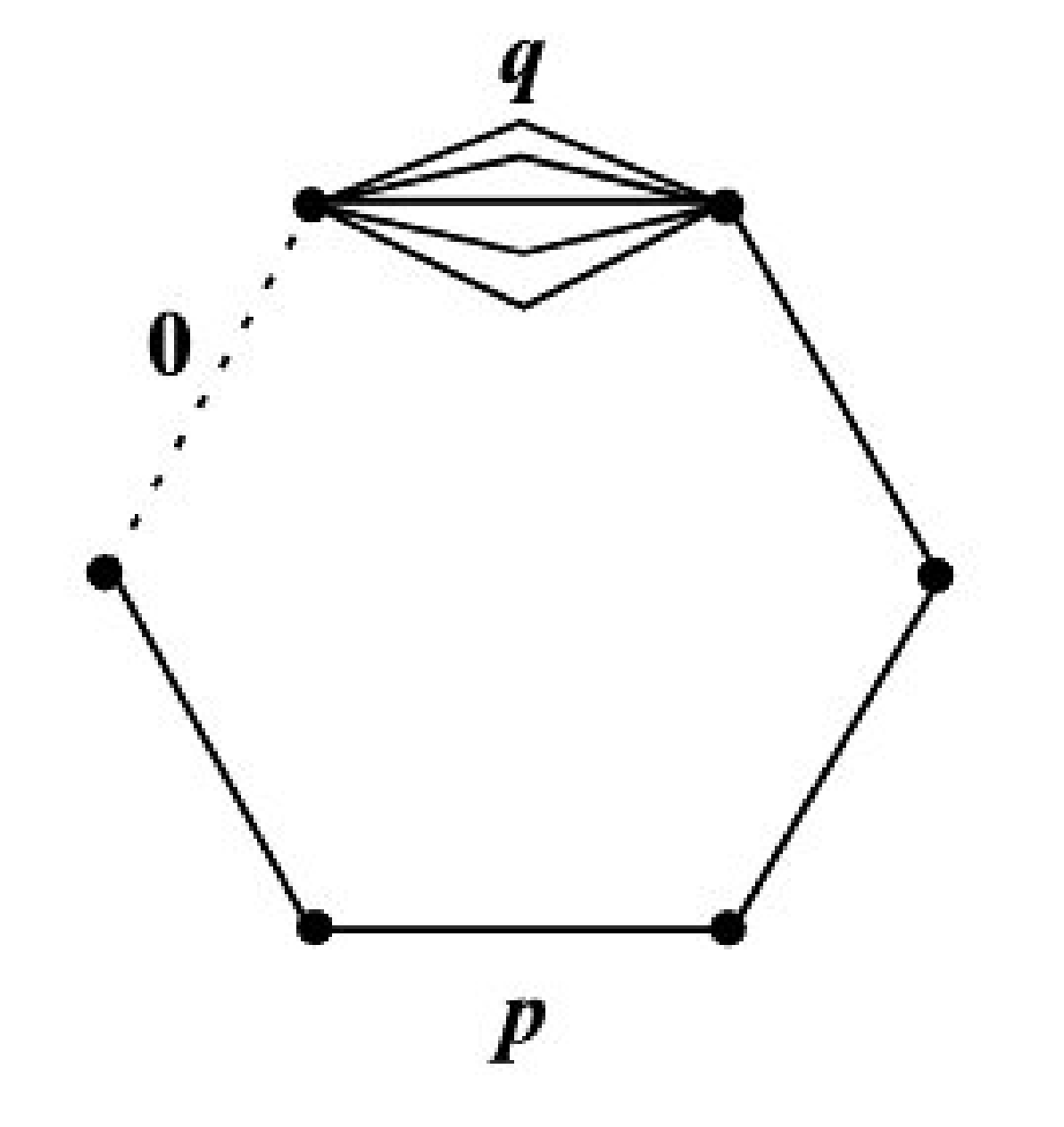,width=1.00in}} \vspace*{8pt}
\caption{Graph $G((i,1^p)\,(1^q))$.\label{f1.3}}
\end{figure}

Next we consider virtual $KL$s of the form $(i,1^p)\,(1^q)$ in the family \emph{\textbf{p\,q}}. The reduced relative Tutte
polynomials of their graphs (Fig. \ref{f1.3}) satisfy the recursion:

$$T(G(i,1^p))=yX^{p-1}T(\overline G(1^q))+xT(G((i,1^{p-1})\,
(1^q)))$$

\noindent where $\overline G(1^q)$ is the dual of the graph
$G(1^q)$, and $T(G((i,1)\,(1)))=xy+x^2y+x^3$, so the general formula
for the reduced relative Tutte polynomial of the graph
$G((i,1^p)\,(1^q))$ is

$$T(G((i,1^p)\,(1^q)))={(x^{q+1}-X^{q+1})(y^{p+1}-Y^{p+1})\over (x-X)(y-Y)} $$
$$-{(y^{p+1}-Y^{p+1}) \over {y-Y}}X^q -{(x^{q+1}-X^{q+1}) \over
{x-X}}Y^n +{(x^q-X^q) \over {x-X}}y^{n+1} +X^qY^n -x^qY^n.$$

\subsection{Family \it{\textbf{p\,1\,q}}}

Based on the link family \emph{\textbf{p\,1\,q}} we construct  three
families of virtual $KL$s with different reduced graphs:
$(i,1^p)\,1\,(1^q)$, $(i,1^p)\,1\,(i,1^q)$, and $(1^p)\,i\,(1^q)$.

\begin{figure}[th]
\centerline{\psfig{file=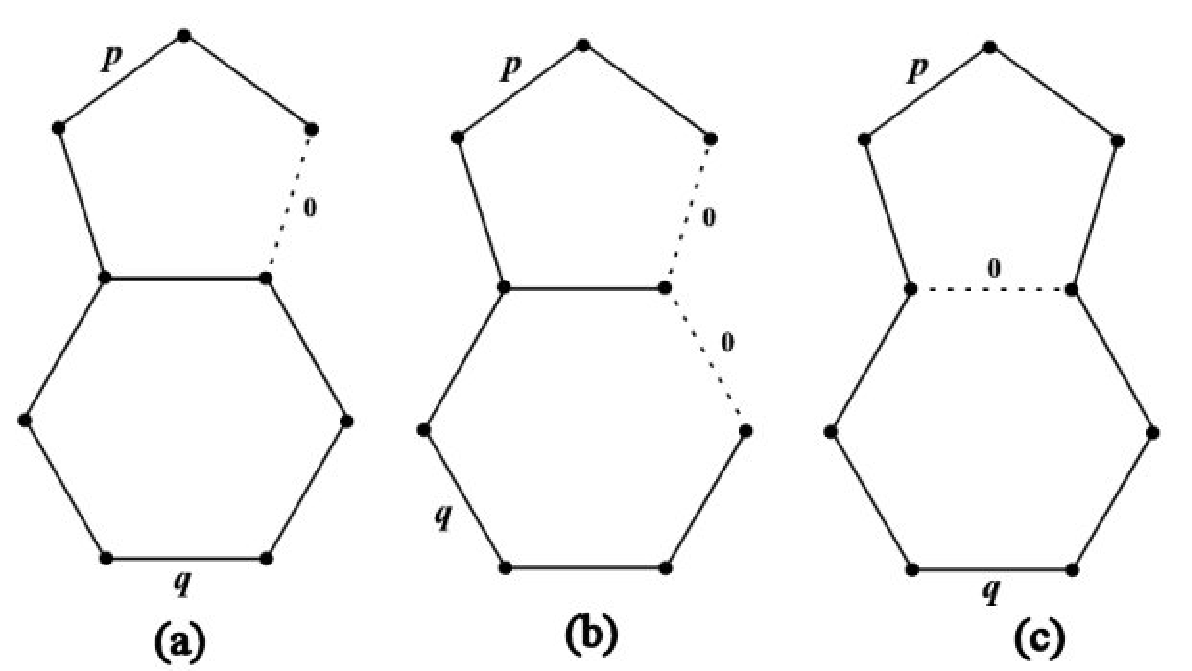,width=2.80in}}
\vspace*{8pt} \caption{Graphs (a) $G((i,1^p)\,1\,(1^q))$; (b) $G((i,1^p)\,1\,(i,1^q))$, and (c) $G((1^p)\,i\,(1^q))$.\label{f1.4}}
\end{figure}

The reduced relative Tutte polynomials of the graphs corresponding
to the family $(i,1^p)\,1\,(1^q)$ (Fig. \ref{f1.4}a) satisfy the
relation:

\begin{equation} \label{for5}
  T(G((i,1^p)\,1\,(1^q)))= yT(G((i,1^{p+q}))+xT(G((i,1^p)))T(G((1^q))),
\end{equation}

\noindent and the reduced relative Tutte polynomials of the graphs
corresponding to the family $(i,1^p)\,1\,(i,1^q)$ (Fig. \ref{f1.4}b)
satisfy the relation:

\begin{equation}\label{for6}
T(G((i,1^p)\,1\,(i,1^q)))= yT(G((1^{p+q}))+xT(G((i,1^p)))T(G((i,1^q))),
\end{equation}

\noindent so the general formulas for their reduced relative Tutte
polynomials can be derived from the previous formulas (\ref{for6}).

The relative Tutte polynomial of the reduced graphs corresponding
to the family $(1^p)\,i\,(1^q)$ (Fig. \ref{f1.4}c) satisfy the
recursion:

\begin{equation}\label{for6red}
T(G((1^p)\,i\,(1^q)))= yX^{q-1}T(G((i,1^p)))+xT(G((1^p)\,i\,(1^{q-1}))),
\end{equation}

\noindent with $T(G((1^p)\,i\,1))=T(G((i,1^p)\,1)$.

Hence, the general formula for the reduced relative Tutte polynomial
of the graph $G((1^p)\,i\,(1^q))$ is

$$T(G((1^p)\,i\,(1^q)))={(x^p-X^p) \over {x-X}}(x+X)x^{q-1}y + {(x^{q-p-1}-X^{q-p-1}) \over
{x-X}}x^pX^{p+1}y$$ $$+{(x^{p-2}-X^{p-2}) \over {x-X}} {(x^q-X^q)
\over {x-X}} x^2y^2+  {(x^q-X^q) \over
{x-X}}(x+X)X^{p-2}y^2+x^{p+q-1}Y.$$

Same as before, from the general formula for reduced relative Tutte
polynomials with negative values of parameters we obtain reduced
relative Tutte polynomials expressed as Laurent polynomials. For
example, for $p=4$, $q=-3$ from the preceding general formula we
obtain reduced relative Tutte polynomial of the virtual knot
$(1,1,1,1)\,i\,(-1,-1,-1)$:

$${-x^3y\over X^3} - {x^2y\over X^2} - {xy\over X} +
{Xy\over x} + {X^2y\over x^2} + {X^3y\over x^3} - 3{y^2\over x} -
{x^2y^2\over X^3} - 2{xy^2\over X^2} - 3 {y^2\over X} - 2 {X
y^2\over x^2} - {X^2y^2\over x^3} + Y$$

\noindent and for $p=-4$, $q=3$ we obtain reduced relative Tutte
polynomial of its mirror image $(-1,-1,-1,-1)\,i\,(1,1,1)$:

$${x^2y\over X^4} - {xy\over X^3} - {y\over X^2} -
    {y\over xX} + {Xy\over x^3} + {X^2y\over x^4} - 2 {y^2\over x^3} -
{xy^2\over X^4} - 2 {y^2\over X^3} - 3 {y^2\over xX^2} - 3 {y^2\over
x^2X} - {Xy^2\over x^4} + {Y\over x^2}.$$

Appropriate substitutions of the variables, yield their
Kauffman bracket and Jones polynomials.

\section{Virtual knots with trivial Jones polynomial and the $Z$-move}

For classical knots, the question whether non-trivial knot with unit Jones polynomial exists is still open.
In the case of virtual knots, it is easy
to make infinitely many non-trivial virtual knots with unit Jones
polynomial \cite {17}. Among 2171 prime virtual knots derived from
classical knots with at most $n=8$ crossings, there are 272 knots (about
12\%) with unit Jones and Kauffman bracket polynomial. The smallest
virtual knot with unit Jones polynomial is $1,1,i,-1,i$ (\cite {17},
Fig. 17), that can be simply generalized to an infinite family of
different non-trivial virtual knots of the form $1^p,i,(-1)^q,i$,
$p-q=1$, with unit Jones polynomial (where $1^p$ denotes a sequence
$1,\ldots,1$ of the length $p$, and $-1^q$ a sequence $-1,\ldots ,-1$
of the length $q$). In the same way, for $p-q=2k+1$, ($k\ge 1$)
all mutually different virtual knots $1^p,i,(-1)^q,i$ will have the
same Jones polynomial as the classical knots $(2k+1)_1$ ($3_1$, $5_1$,
$7_1$, $\ldots$, or $3$, $5$, $7$, $\ldots $ in Conway notation),
respectively.

For virtual $KL$s that can be reduced to unknot (unlink) by a series
of Reidemeister moves for virtual knots and $Z$-moves we will say
that they are $Z$-move equivalent to the unknot (unlink). The most
of mentioned 272 knots with unit Jones polynomial are $Z$-move
equivalent to the unknot. Hence, R. Fenn, L.H. Kauffman, and V.O. Manturov \cite {VP} proposed $Z$-move
conjecture:

\medskip

{\bf Conjecture 3.1}
{\it Every knot with unit Jones polynomial is $Z$-move equivalent to the
unknot.}

\medskip

\begin{figure}[th]
\centerline{\psfig{file=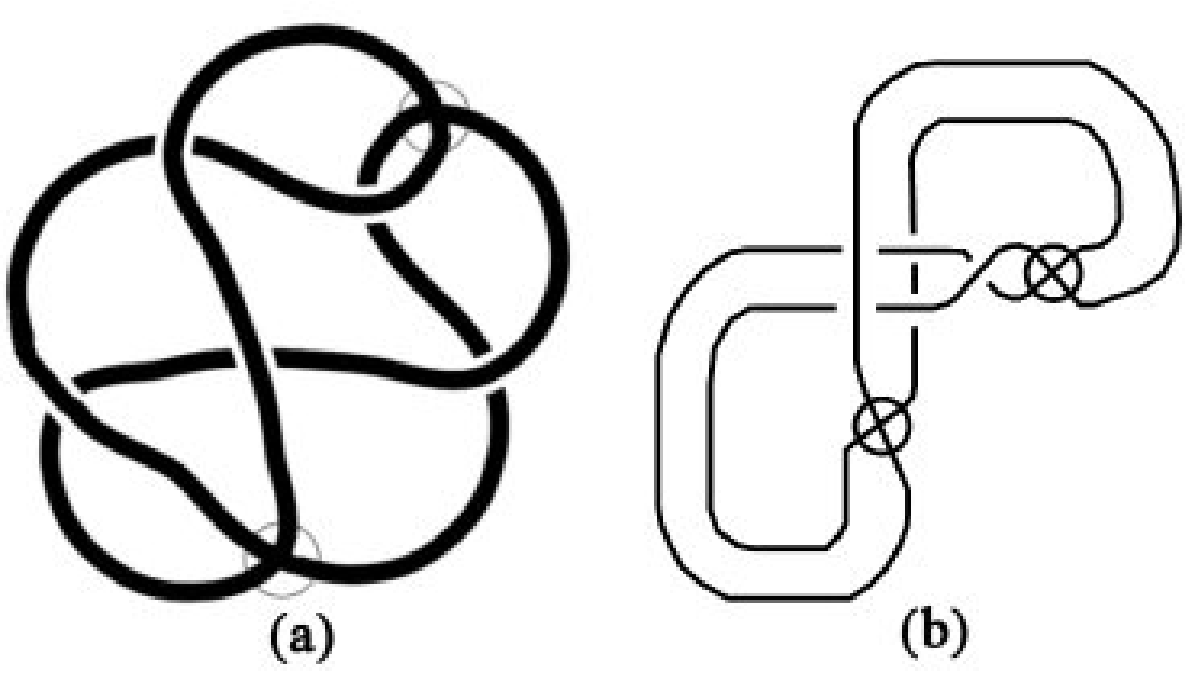,width=2.00in}} \vspace*{8pt}
\caption{(a),(b) Virtual knot
$KS=(((1,(i,1),-1),-1),i,1)$.\label{f1.7}}
\end{figure}

\begin{figure}[th]
\centerline{\psfig{file=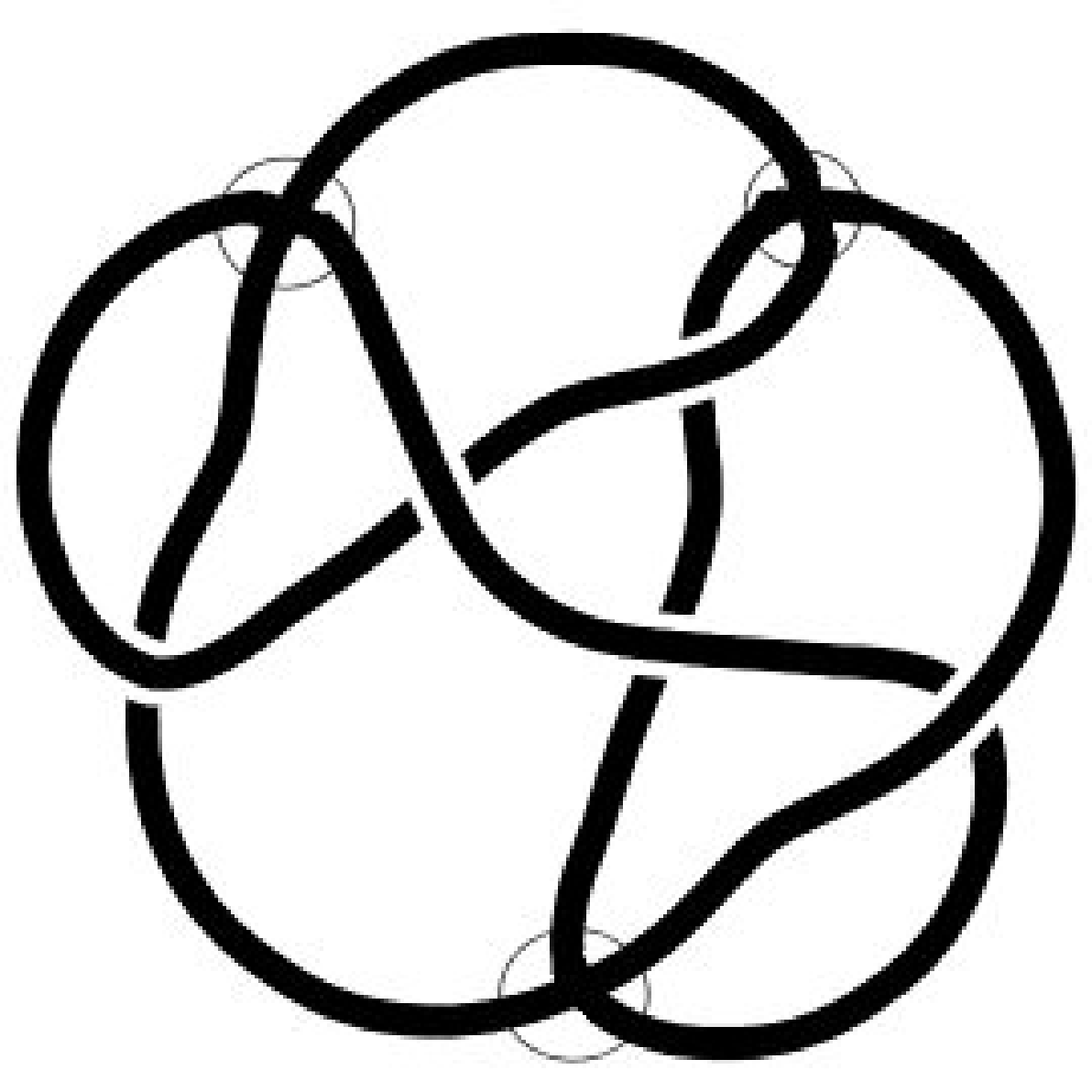,width=1.00in}} \vspace*{8pt}
\caption{Virtual knot $S7=(i,1)\,1,(i,-1)\,-1,(i,1)$.\label{f1.8}}
\end{figure}

The main candidates for counterexamples to this conjecture have been
virtual knots $KS=(((1,(i,1),-1),-1),i,1)$ (Fig. \ref{f1.7}a,b) and
$S7=(i,1)\,1,(i,-1)\,-1,(i,1)$ (Fig. \ref{f1.8}). For each of them,
non-triviality can be proved in many different ways (e.g., for the
first of them by parity arguments, or by computing their other
polynomial invariants: Sawollek polynomial, Miyazawa polynomials
\cite {19}, Dye-Kauffman arrow polynomial \cite {18}, or 2-cabled Jones
polynomial. It is interesting to notice that the knot $KS$ has unit
Miyazawa polynomials as well.

In the next section we prove that the virtual knot $KS$ with unit Jones polynomial is not $Z$-move
equivalent to unknot, giving a counterexample to this conjecture (already  known to be false, see \cite{MP,IVKT}). Methods for proving that counterexamples that we discuss below are indeed counterexamples depend on the use of parity
in virtual knot theory as introduced in \cite{MP,SL,IVKT}. See the next section for a detailed discussion of this method.
\bigbreak

\begin{figure}[th]
\centerline{\psfig{file=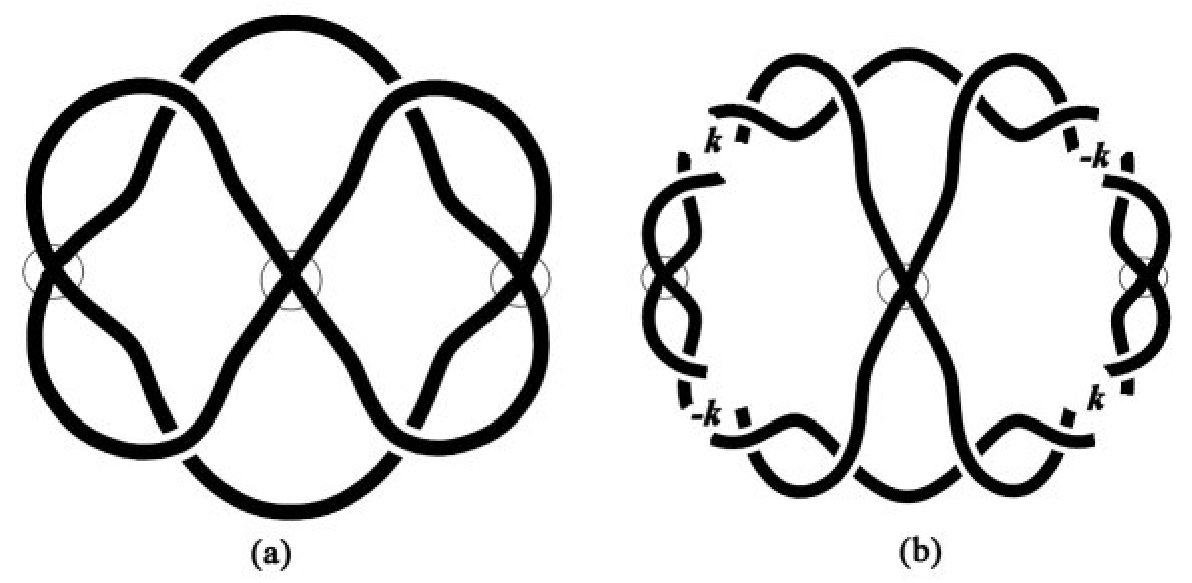,width=2.60in}} \vspace*{8pt}
\caption{(a) Virtual knot $(1,i,-1)\,i\,(1,i,-1)$; (b) family
$(1^k,i,(-1)^k)\,i\,$ $(1^k,i,(-1)^k)$.\label{f1.9}}
\end{figure}

Among virtual knots with unit Jones polynomial, probably the most
interesting is the virtual knot $(1,i,-1)\,i\,(1,i,-1)$ (Fig.
\ref{f1.9}a), which generates the family of virtual knots
$(1^k,i,(-1)^k)\,i\,(1^k,i,(-1)^k)$ ($k\ge 1$) (Fig. \ref{f1.9}b)
which cannot be distinguished from unknot by any of the mentioned
polynomial invariants, except 2-colored Jones polynomial. All
members of this family are $Z$-move equivalent to the unknot.

\begin{figure}[th]
\centerline{\psfig{file=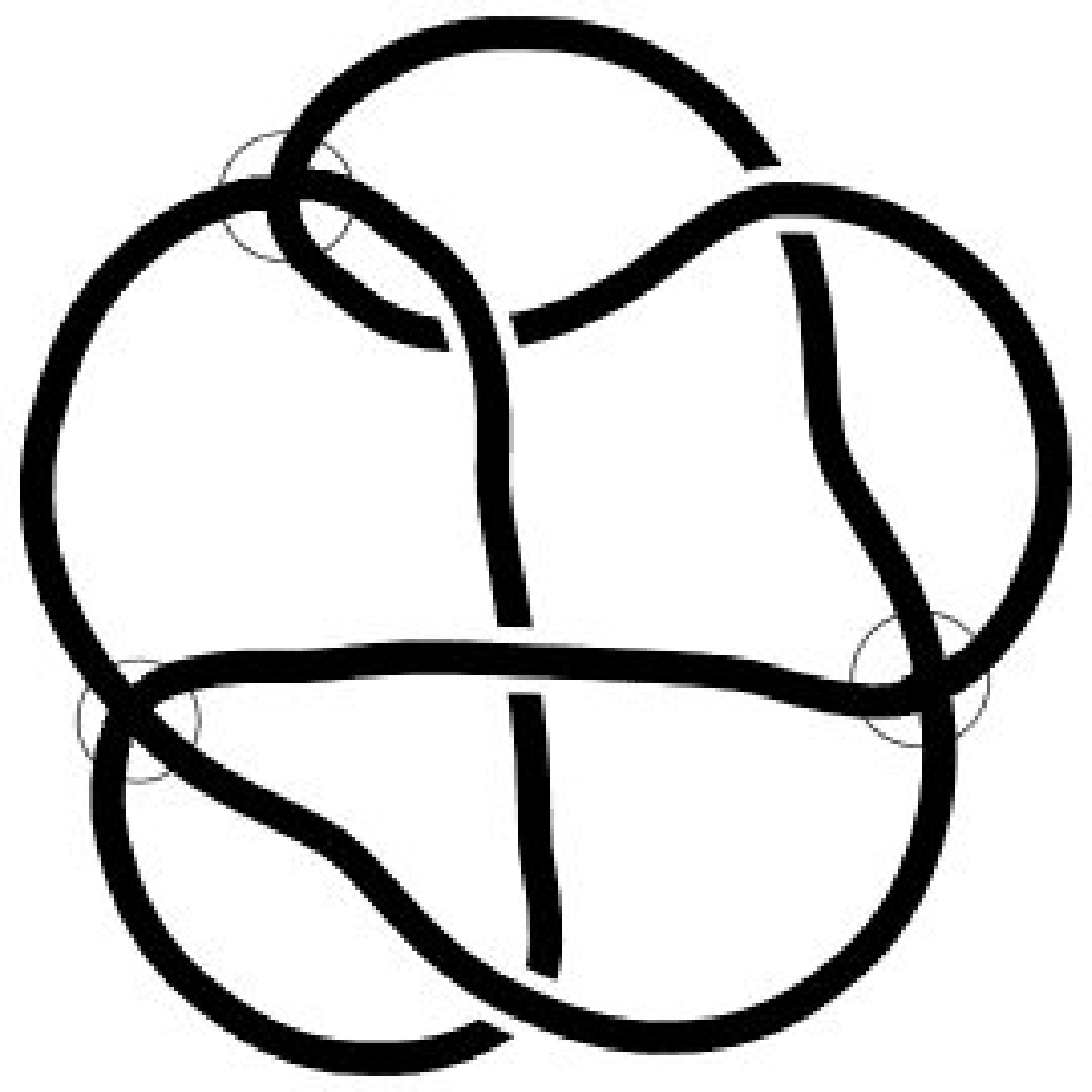,width=1.30in}} \vspace*{8pt}
\caption{Virtual knot $(i,1)\,(1,i)\,1\,(-1,i)$.\label{f1.10}}
\end{figure}

Another interesting virtual knot is $(i,1)\,(1,i)\,1\,(-1,i)$ which
has all trivial mentioned polynomial invariants, except 2-colored
Jones polynomial. Moreover, it is another candidate for a counterexample
to $Z$-move conjecture (Fig. \ref{f1.10}).

Using $Z$-move reduction we can prove many simple facts about Jones
polynomial of virtual $KL$s. For example:

\begin{itemize}
\item if every twist of a knot or link $L$ which contains virtual
crossings has an even number of them, the Jones polynomial of $L$ is
equal to the Jones polynomial of the classical link $L'$ obtained from
$L$ by deleting virtual crossings.

\item Let be given an alternating knot $K$ with $n$ crossings. In
order to unknot it, in any minimal diagram it is sufficient to
make at most $m=[{n \over 2}]$ crossing changes. If a minimal
diagram $D$ of $K$ can be unknotted by the crossing changes in the
crossings $C_1$, $\ldots $, $C_k$ ($k\le m$), the virtual knot
diagram $D'$, obtained from $D$ by substituting every positive
crossing $C_i$ ($i\in \{1,\ldots ,k\}$) by $i,-1,i$, and every
negative crossing $C_j$ ($j\in \{1,\ldots ,k\}$) by $i,1,i$, has
unit Jones polynomial. Are all virtual knots obtained in this way
non-trivial? Can we obtain from two different minimal diagrams of
$K$ two different virtual knots?
\end{itemize}

\section{The Parity Bracket Polynomial}

In this section we introduce the Parity Bracket Polynomial of Vassily Manturov \cite{MP}.
This is a generalization of the bracket polynomial to virtual knots and links that uses the parity of the crossings. A crossing is {\it odd}  (of odd parity) if
it flanks an odd number of symbols in the Gauss code of that diagram. A crossing that flanks an even number of symbols is said to be {\em even.} For example, if we have a Gauss code of the form $1212$ (here we have specified only the crossings $1$ and $2$, and nothing about
their signs or whether they occur over or under) then both crossings are odd. On the other hand in the code $123123$, all crossings are even and in the code $123213$, the crossings $1$ and $2$ are odd, while the crossing $3$ is even. In a classical knot diagram, it is easy to see (by the Jordan curve theorem) that all crossings are even. But in a virtual knot diagram it is possible to have odd crossings.

\begin{figure}[th]
\centerline{\psfig{file=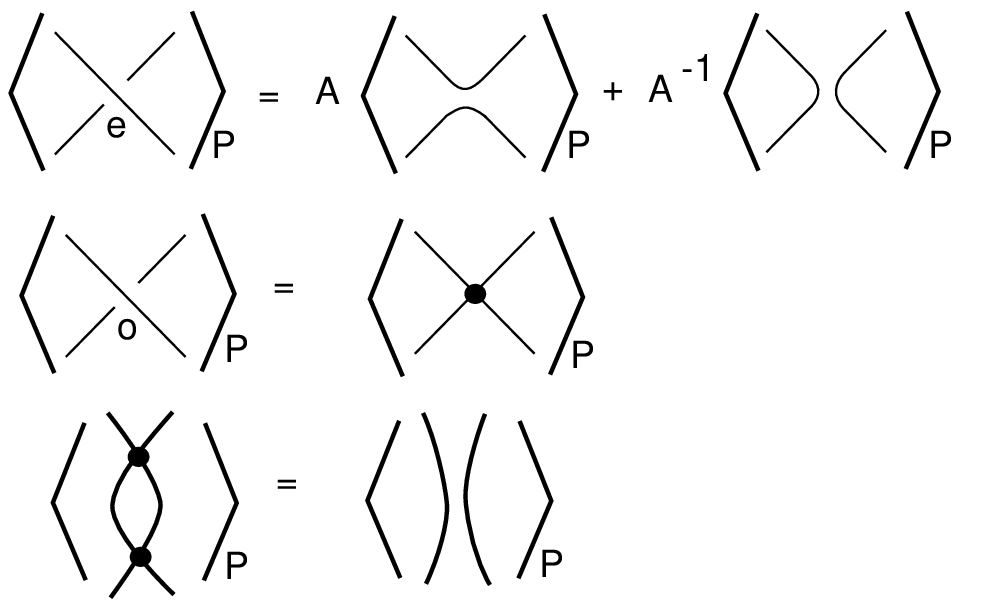,width=2.20in}} \vspace*{8pt}
\caption{Parity bracket expansion.\label{f1.11}}
\end{figure}

We define a {\em Parity State} of a virtual diagram $K$ to be a labeled virtual
graph obtained from $K$ as follows: For each odd crossing in $K$ replace the crossing by a graphical node. For each even crossing in $K$ replace the crossing by one of its two possible smoothings, and label the smoothing site by $A$ or $A^{-1}$ in the usual way. Then we define the parity bracket by the
state expansion formula
$$\langle K \rangle _{P} = \sum_{S}A^{n(S)}[S]$$
where $n(S)$ denotes the number of $A$-smoothings minus the number of $A^{-1}$ smoothings and
$[S]$ denotes a combinatorial evaluation of the state defined as follows: {\em First reduce the state
by Reidemeister two moves on nodes as shown in Fig.\ref{f1.11}.} The reader should note that we have labeled even crossings by {\em e} and odd crossings by {\em o}. There should be no confusion between this notation and the notation we used previously for a virtual edge in a knot graph.

\begin{figure}[th]
\centerline{\psfig{file=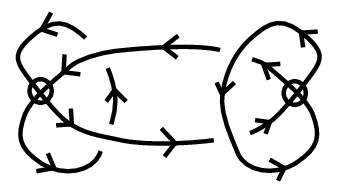,width=2.00in}} \vspace*{8pt}
\caption{Kishino diagram.\label{f1.12}}
\end{figure}

Here the graphs are taken up to virtual equivalence (planar isotopy plus detour moves on the virtual crossings. See \cite{IVKT}.). We regard the reduced state
as a disjoint union of standard state loops (without nodes) and graphs that irreducibly contain nodes.
With this we write $$[S] = (- A^{2} - A^{-2})^{l(S)} [G(S)]$$ where $l(S)$ is the number of standard loops in the reduction of the state $S$ and $[G(S)]$ is the disjoint union of reduced graphs that contain nodes.
In this way, we obtain a sum of Laurent polynomials in $A$ multiplying reduced graphs as the Manturov Parity Bracket. It is not hard to see that this bracket is invariant under regular isotopy and detour moves
and that it behaves just like the usual bracket under the first Reidemeister move. However, the use of
parity to make this bracket expand to graphical states gives it considerable extra power in some situations. For example, consider the Kishino diagram in Fig.\ref{f1.12}. We see that all the classical crossings in this knot are odd. Thus the parity bracket is just the graph obtained by putting nodes at each of these crossings. The resulting graph does not reduce under the graphical Reidemeister two moves, and so we conclude that the Kishino knot is non-trivial and non-classical. Since we can apply the parity
bracket to a flat knot by taking $A = -1$, we see that this method shows that the Kishino flat is non-trivial.

\bigskip
\noindent {\bf Remark.} It should be mentioned that the graph-link theory and free knot theory  of Manturov-Ilyutko is the best setting for the Manturov bracket since it  is abstract, graphical and does not feel the Z-move since the free knot theory is already invariant under the $Z$-move.
The interested reader should see \cite{MP} for more details about this theory. Here, we do not use the free knot theory but rather, we use standard virtual knot theory and we formulate a version of the parity bracket in this context. This formulation is given in  \cite{IVKT}.  The basic idea of the parity bracket is
due to Manturov.

\begin{figure}[th]
\centerline{\psfig{file=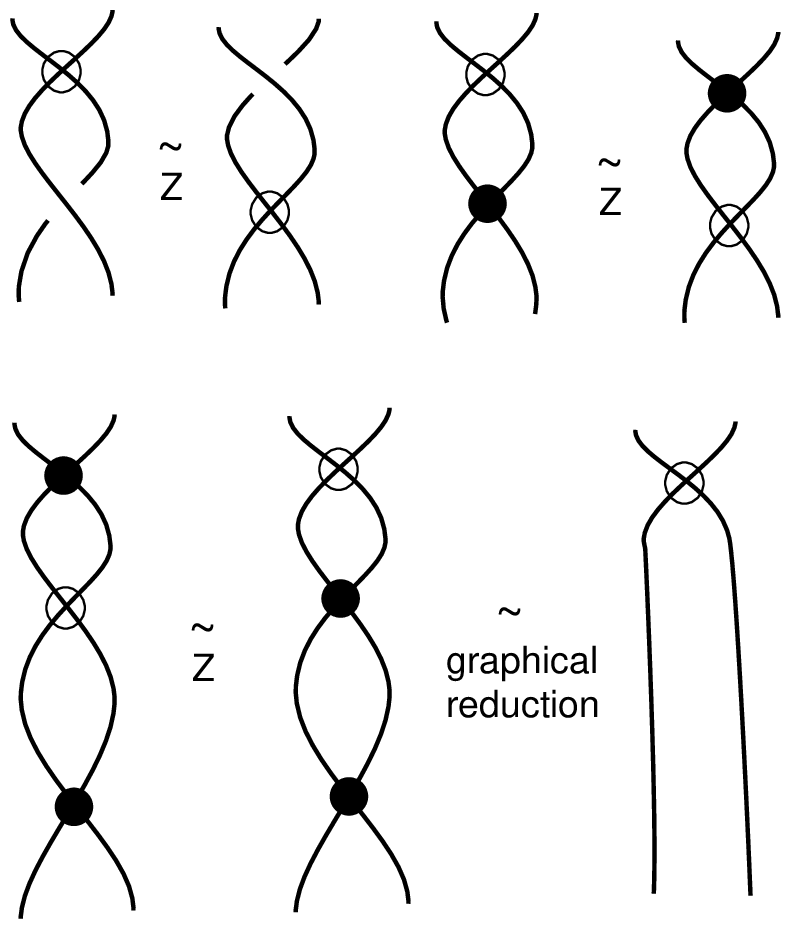,width=2.00in}} \vspace*{8pt}
\caption{$Z$-move and graphical $Z$-move.\label{f1.13}}
\end{figure}

\begin{figure}[th]
\centerline{\psfig{file=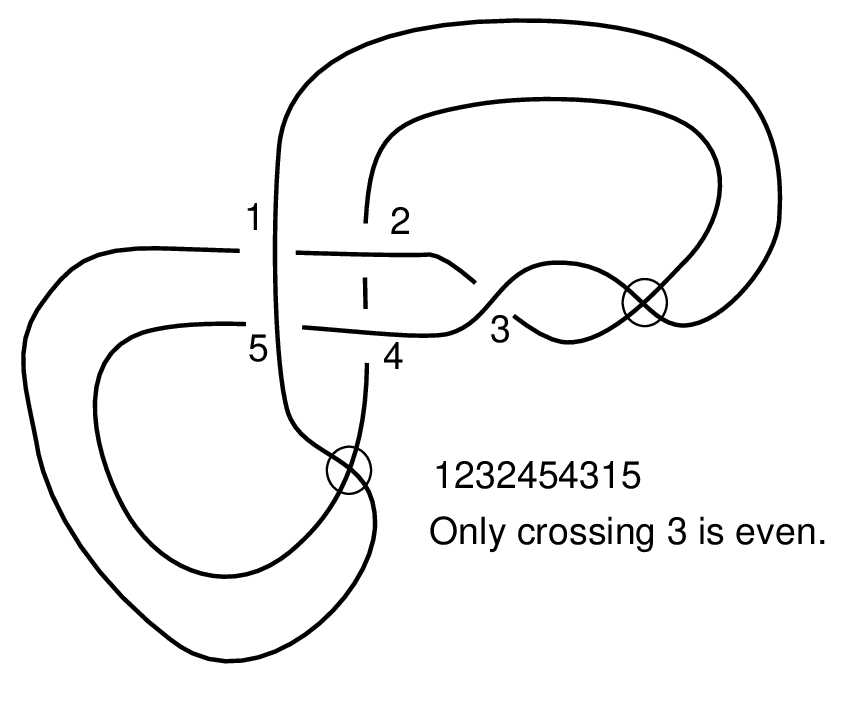,width=1.80in}} \vspace*{8pt}
\caption{A knot $KS$ with unit Jones polynomial.\label{f1.14}}
\end{figure}

\begin{figure}[th]
\centerline{\psfig{file=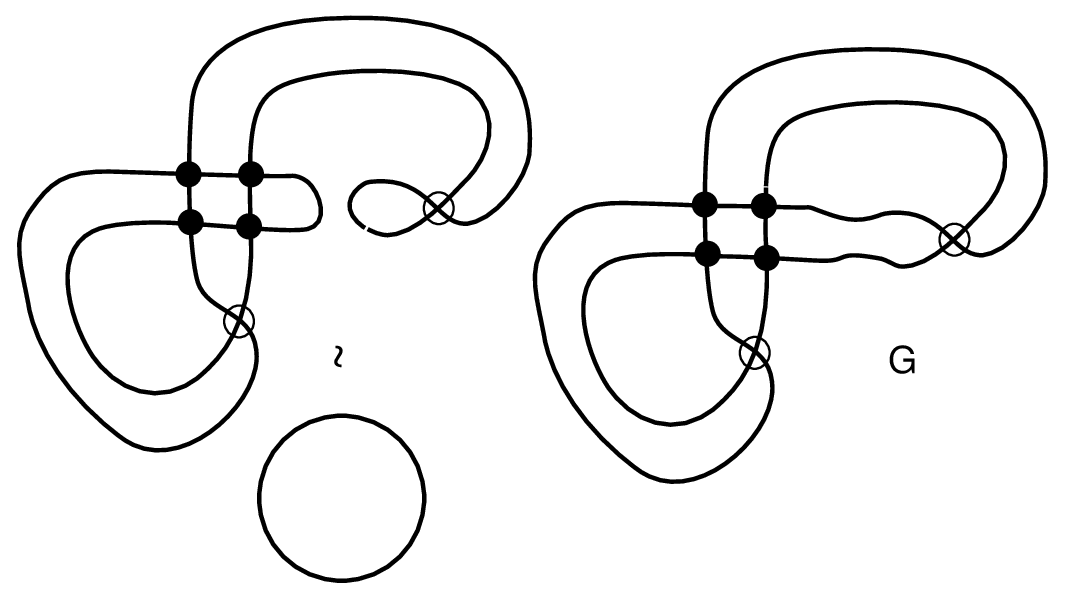,width=2.20in}} \vspace*{8pt}
\caption{Parity bracket states for the knot $KS$.\label{f1.15}}
\end{figure}

In  Fig. \ref{f1.13} we illustrate the {\it $Z$-move}  and the {\it graphical $Z$-move}. Two virtual knots or links
that are related by a $Z$-move have the same standard bracket polynomial. This follows directly from our discussion in the previous section. We would like to analyze the structure of $Z$-moves using the parity
bracket. In order to do this we need a version of the parity bracket that is invariant under the $Z$-move.
In order to accomplish this, we need to add a corresponding $Z$-move in the graphical reduction process for the parity bracket. This extra graphical reduction is indicated in Fig. \ref{f1.13} where we show a graphical $Z$-move. The reader will note that graphs that are irreducible without the graphical $Z$-move can become reducible if we allow graphical $Z$-moves in the reduction process. See Fig.\ref{f1.13} for an example of this process as well as an illustration of the graphical $Z$-move. For example, the graph associated with the Kishino knot is reducible under graphical $Z$-moves. However, there are examples of
graphs that are not reducible under graphical $Z$-moves and Reidemeister two moves. An example of such a graph occurs in the parity bracket of the knot $KS$ shown in  Fig. \ref{f1.14}. This knot has one even classical crossing and four odd crossings. One smoothing of the even crossing yields a state that reduces to a loop with no graphical nodes, while the other smoothing yields a state that is irreducible even when the $Z$-move is allowed (see  Fig. \ref{f1.15}). The upshot is that this knot KS is not $Z$-equivalent to any classical knot. Since one can verify that $KS$ has unit Jones polynomial, this example is a counterexample to a conjecture of Fenn, Kauffman and Manturov \cite{VP} that suggested that a knot with unit Jones polynomial should be $Z$-equivalent to a classical knot. The existence of such counterexamples via parity was first pointed out by Vassily Manturov in 2009.

Parity is clearly an important theme in virtual knot theory and will figure in many future investigations of this subject. The type of construction that we have indicated for the bracket polynomial in this section
can be varied and applied to other invariants. Furthermore the notion of describing a parity for crossings
in a diagram is also susceptible to generalization. For more on this theme the reader should consult
\cite{MP1,PT} and \cite{SL} for our original use of parity for another variant of the bracket polynomial.

\section{Portraits of families of virtual $KL$s}

Recursive and general formulas for the reduced relative Tutte
polynomials can be computed for different families of virtual $KL$s,
and from them we obtain Jones polynomials and Kauffman bracket
polynomials of the considered families of virtual $KL$s.

\begin{figure}[th]
\centerline{\psfig{file=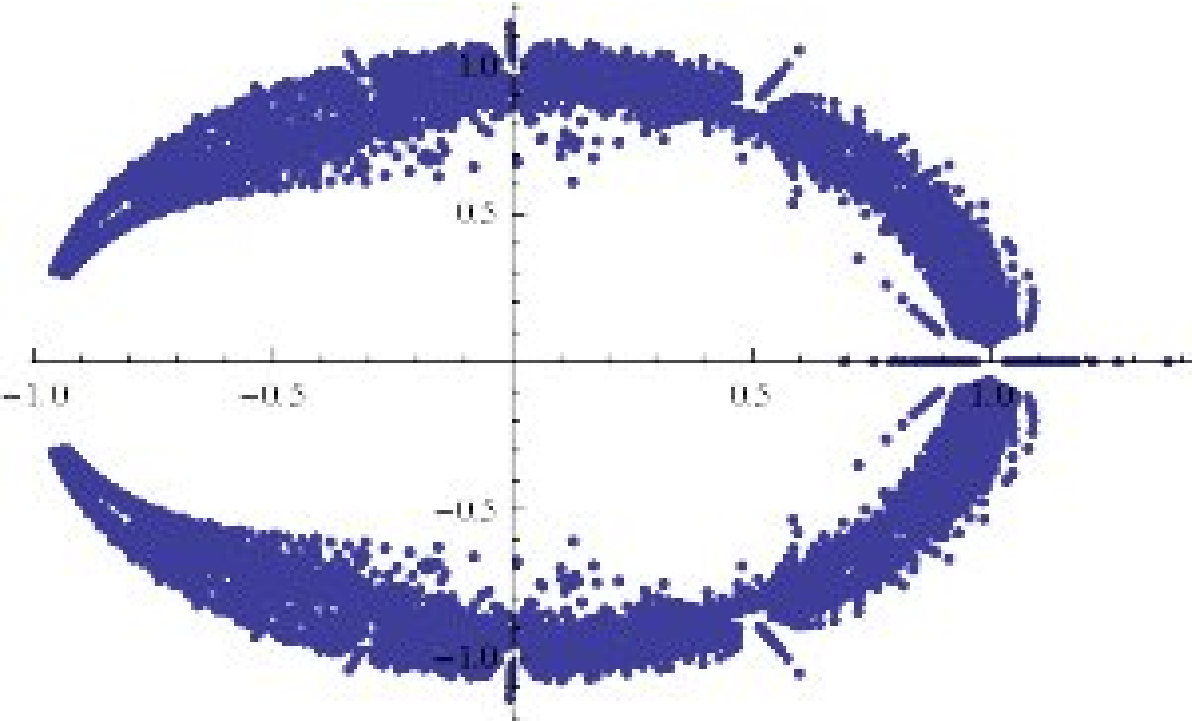,width=2.10in}} \vspace*{8pt}
\caption{Plot of the zeroes of Jones polynomial for virtual $KL$
family $(i,1^p)\,(1^q)$ ($1\le p\le 20$, $2\le q\le
20$).\label{f1.16}}
\end{figure}

\begin{figure}[th]
\centerline{\psfig{file=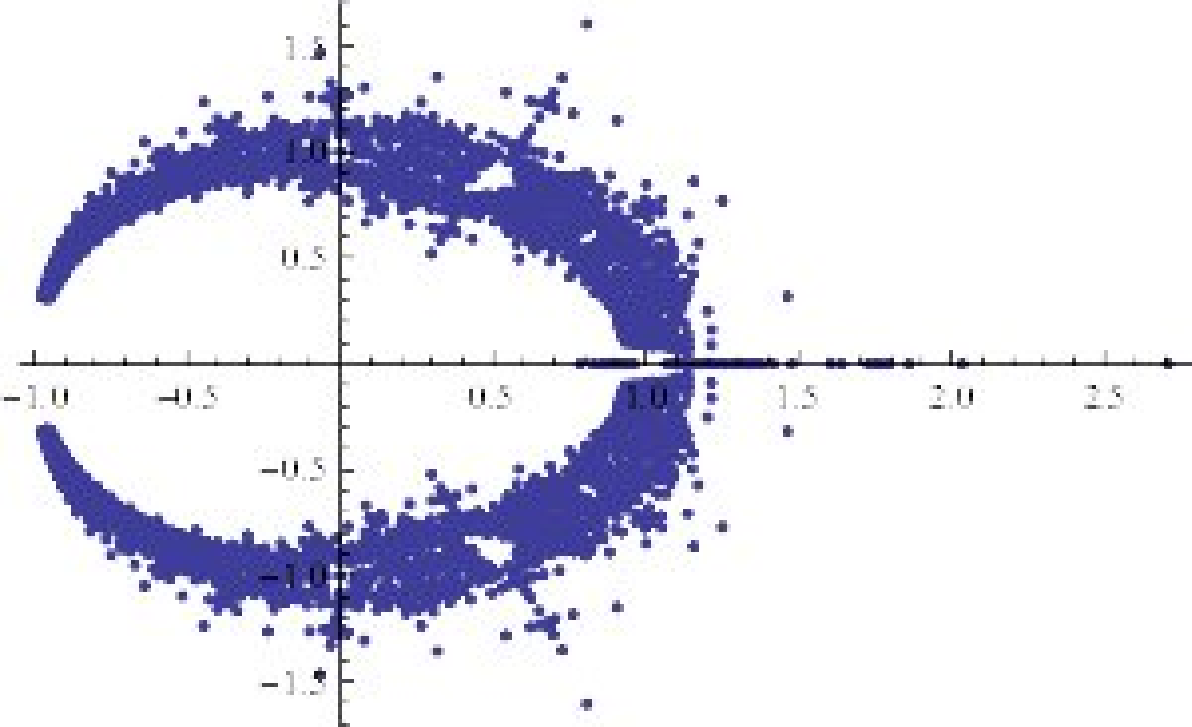,width=3.20in}} \vspace*{8pt}
\caption{Plot of the zeroes of Jones polynomial for virtual $KL$
family $(1^p)\,i\,(1^q)$ ($2\le p\le 20$, $2\le q\le
20$).\label{f1.17}}
\end{figure}

\begin{figure}[th]
\centerline{\psfig{file=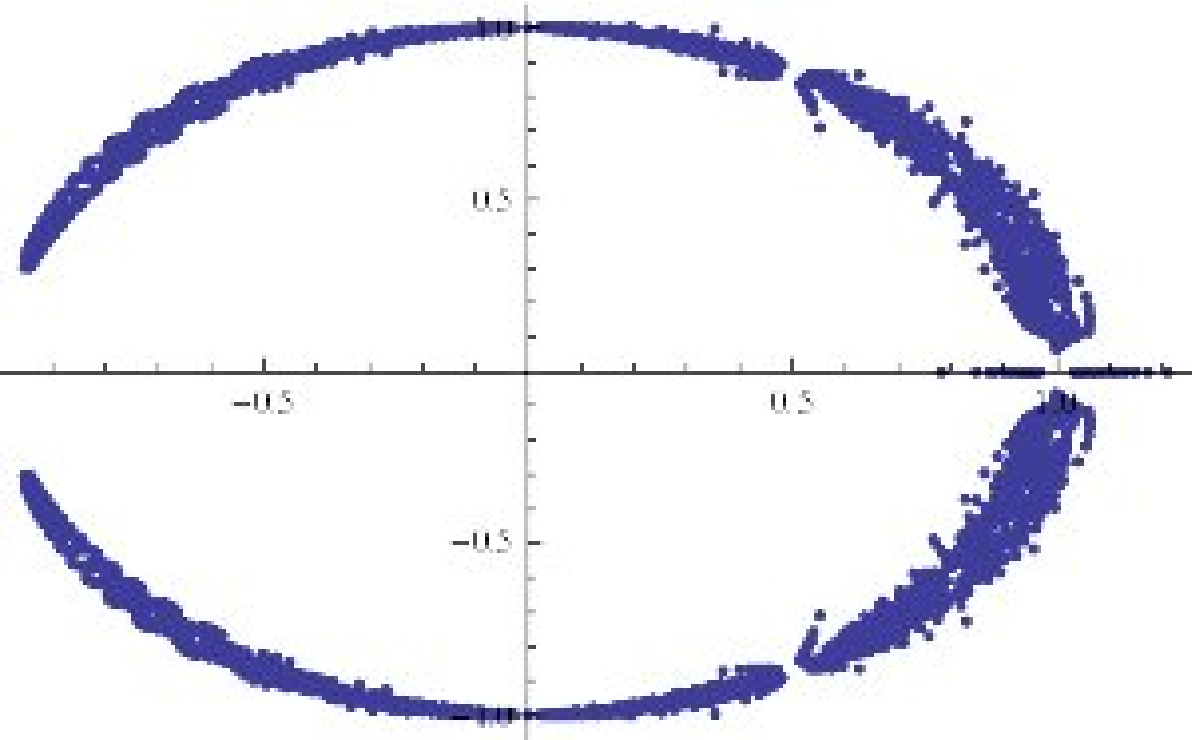,width=2.00in}} \vspace*{8pt}
\caption{Plot of the zeroes of Jones polynomial for virtual $KL$
family $(1^p)\,i\,((-1)^q)$ ($2\le p\le 20$, $2\le q\le
20$).\label{f1.18}}
\end{figure}

Obtained results can be used to study properties of the reduced
relative Tutte polynomials of virtual $KL$s and zeros of Kauffman
bracket polynomials and Jones polynomials. The plot of zeros of
Jones polynomials of virtual $KL$ family is specific to the family
and will be referred to as the ``portrait of a virtual link family''.

The portrait of the virtual link family
$(i,1^p)\,(1^q)$ ($1\le p\le 20$, $2\le q\le 20$) is shown in Fig.
\ref{f1.16}.\\

Portraits of families of virtual $KL$s obtained for different
choices of signs of parameters are compared in Figs. 15 and 16. Fig.
\ref{f1.17} is the portrait of the family $(1^p)\,i\,(1^q)$ for
$2\le p\le 20$, $2\le p\le 20$, and the Fig. \ref{f1.18} corresponds
to the family $(1^p)\,i\,((-1)^q)$ for  $p,q \in Z^+$. More detailed
results of this kind will be given in the forthcoming paper.

\end{document}